\documentclass{amsart}

\usepackage[all]{xy}
\usepackage{amssymb,latexsym}
\usepackage{amsmath}
\usepackage{array}
\usepackage{eqnarray}
\usepackage[utf8]{inputenc}
\usepackage{setspace}
\usepackage{amsthm}
\usepackage{url}
\usepackage{appendix}
\usepackage[margin=1in]{geometry}
\usepackage{todonotes}
\usepackage{amsthm}
\usepackage[parfill]{parskip}
\usepackage{tikz}
\usepackage{tikz-cd}
\usepackage{graphicx}
\graphicspath{ {./images/} }
\makeatletter
\newcommand{\xRightarrow}[2][]{\ext@arrow 0359\Rightarrowfill@{#1}{#2}}
\makeatother

\newcommand{\Sp}{\operatorname{Sp}}

\newcommand{\comment}[1]{}

\newtheorem{defn0}{Definition}[section]
\newtheorem{prop0}[defn0]{Proposition}
\newtheorem{thm0}[defn0]{Theorem}
\newtheorem{lemma0}[defn0]{Lemma}
\newtheorem{corollary0}[defn0]{Corollary}
\newtheorem{example0}[defn0]{Example}
\newtheorem{remark0}[defn0]{Remark}
\newtheorem{conjecture0}[defn0]{Conjecture}
\newtheorem{ques0}[defn0]{Question}
\newtheorem*{mthm0}{Main Theorem}

\newenvironment{proposition}{\bigskip \begin{prop0}}{\end{prop0}}

\newenvironment{corollary}{\bigskip \begin{corollary0}}{\end{corollary0}}

\newcommand{\ii}{\textbf{i}}
\newcommand{\jj}{\textbf{j}}
\newcommand{\kk}{\textbf{k}}
\newcommand{\Z}{\mathbb{Z}}
\newcommand{\R}{\mathbb{R}}
\newcommand{\HH}{\mathbb{H}}
\newcommand{\OO}{\mathbb{O}}

\newcommand{\C}{\mathbb{C}}
\newcommand{\CP}{\mathbb{C}P}

\renewcommand{\ker}{\operatorname{Ker}}

\newcommand{\Ogroup}{\operatorname{O}}
\newcommand{\SO}{\operatorname{SO}}
\newcommand{\Spin}{\operatorname{Spin}}
\newcommand{\U}{\operatorname{U}}
\newcommand{\SU}{\operatorname{SU}}
\newcommand{\Gl}{\operatorname{GL}}
\newcommand{\eG}{\operatorname{G}_2}

\newcommand{\Ad}{\operatorname{Ad}}

\newcommand{\Stab}{\operatorname{Stab}}

\usepackage{soul}

\begin{document}
\title{G-invariant Spin Structures on Spheres}

\author{Jordi Daura Serrano}
\address{Jordi Daura Serrano, Department de Màtematiques i Informàtica, Universitat de Barcelona (UB), Gran Via de les Corts Catalanes 585, 08007 Barcelona (Spain)}
\email{jordi.daura@ub.edu}
\author{Michael Kohn} 
\address{Michael Kohn, Centre of Mathematics, Wilberforce rd, CB3WA Cambridge (UK)}
\email{michael.kohn18@imperial.ac.uk}
\author{Marie-Amélie Lawn}
\address{Marie-Amélie Lawn, Department of Mathematics, Imperial College, 180 Queen's Gate, London (UK) }
\email{mlawn@ic.ac.uk}

\dedicatory{This paper is dedicated to Lionel B\'erard Bergery, who was one of the greatest experts in homogeneous spaces, and whose enthusiasm and ability to communicate his love to the subject remains unmatched.}  
\pagenumbering{arabic}
\begin{abstract}
We examine which of the compact connected Lie groups that act transitively on spheres of different dimensions leave the unique spin structure of the sphere invariant. We study the notion of invariance of a spin structure and prove this classification in two different ways; through examining the differential of the actions and through representation theory.  
\end{abstract}
\date{}
\maketitle
\section*{Introduction}

Given a Lie group $G$ acting on a manifold $M$, it is a natural question to ask which structures of the manifold are preserved by $G$. For example, if $M$ is orientable, connected Lie groups always preserve an orientation. However, if we consider that $M$ admits spin structures (which can be regarded as a refinement of orientation) the question of their preservation by $G$ is more complicated. We will show in the main theorem that even if the manifold has a unique spin structure and the Lie group is connected, it is possible for the action to not preserve it. 

The question is even more relevant in the case of homogeneous spaces $G/H$, where the group action determines the manifold. Here, spin structures can be nicely characterized in terms of lifts of the isotropy representation to the group $\Spin(n)$ (see \cite{ALEKSEEVSKYD} and references therein). 

The question is motivated by the utility of $G$-invariant spin structures; in addition to their elegant form mentioned above, G-invariance is crucial to compute the spinors which are invariant with respect to the connection induced by the action of the Lie group. This fact is a further motivation to carefully study the case of homogeneous spaces: Most examples of geometries admitting special spinors are indeed homogeneous (see for example \cite{wang1989parallel}, \cite{bar1993real} and \cite{Agr}).   

A particularly interesting case are the spheres. Most commonly viewed as the homogeneous space  $S^n=\SO(n+1)/\SO(n)$, they can actually be realized, due to their great amount of symmetries, as various homogeneous decompositions, according to the different groups acting transitively and effectively on them. These groups were classified by D.Montgomery and H.Samelson (see \cite{Samuelson}). 
Part of the preliminaries and an appendix at the end of this paper are dedicated to describing each of these nine possible group actions, together with their isotropy representations, as a comprehensible survey of these classical results is very difficult to find in the literature. Understanding the features of these actions is of particular importance in differential geometry; one of their most remarkable properties being the well-known fact that the transitive effective Lie  groups on spheres appear as holonomy groups according to Berger's classification (\cite{Berger}) on simply connected Riemannian manifolds, with the exception of $\Spin(9)$ and $\Sp(n)\cdot \U(1)$.

Our main theorem studies which of these actions leave the unique spin structure of spheres invariant. More precisely, we prove the following:

\begin{mthm0}
	Let $G$ be a compact connected Lie group acting transitively and effectively on a sphere of appropriate dimension. Then, its unique spin structure is $G$-invariant if and only if $G$ is simply connected, or for $n$ odd,  $\Sp(n+1)\cdot\U(1)$ or $\Sp(n+1)\cdot\Sp(1)$.
\end{mthm0}

Perhaps the most surprising fact is that there are non-simply connected Lie groups which preserve the spin structure. 

We present two different methods to prove the main theorem. One is of more differential geometric nature and based on the definition of the isotropy representation using differentials, through the fact that the actions are linear. The other approach uses the characterization of the isotropy representation as a restriction of the adjoint representation of the  transitive Lie group, which enables the use of the tools of representation theory. This method is more general as it is not necessary for the action to be linear. It is important to point out that in general we do not need to compute the whole isotropy representation, since it is enough to find the image of loops whose classes generate the fundamental group of the stabiliser $H$. By choosing the adequate loop representative we can reduce the complexity of the computation. 

The paper is divided as follows. First, we recall the basic concepts about the theory of  homogeneous spaces we need, focusing on describing the different transitive group actions on spheres, and recall the concept of a $G$-invariant spin structure, with special attention to the case of homogeneous spaces. In the second section, we prove our main result, which we divide into lemmata for each of the transitive group actions. We end the paper with a few closing remarks and results. In particular we point out that our main theorem yields a complete classification of $G$-invariant spin structures on spheres, with $G$ a compact connected Lie group acting transitively.

\subsection*{Acknowledgements} The authors would like to thank Travis Schedler for fruitful discussions about the representation theory used in this paper and Ilka Agricola and Dmitri Alekseevsky for helpful comments, as well as the anonymous referee for his suggestions.

\subsection*{Data Availability Statement} No datasets were generated or analysed during the current study.

\section{Preliminaries}
\subsection{Homogeneous spaces}

In this section, we recall some well-known facts about homogeneous spaces. For more details we refer  to \cite[Chapter 7]{BesseArthur2008EM}. Let $M$ be an orientable $n$-dimensional Riemannian manifold on which a connected Lie group $G$ acts transitively from the left (i.e. a homogeneous space) and fix a point $o\in M$. Then, it is well-known that $M \simeq G/H$, where $H=\Stab(o)$. We denote by $\mu_{g'}: M\rightarrow M,\,gH\mapsto g'gH$ the group action of $G$ on $M$. Additionally, we assume that the isotropy subgroup $H$ is compact, where $o:=eH$. Then $M$ admits an invariant Riemannian metric, and the group $G$ acts by orientation-preserving isometries. 

 
We denote by $\Ad^G: G \rightarrow  \Gl(\mathfrak{g})$ and $\Ad^H: H \rightarrow  GL(\mathfrak{h})$ the adjoint representations of $G$ and $H$ respectively,
 where $\mathfrak{g}=T_eG$ denotes as usual the Lie algebra of $G$ and $\mathfrak{h}$ the Lie algebra of $H$. Riemannian homogeneous spaces are reductive, i.e. there exists an $H$-invariant vector space $\mathfrak{m}$ such that $\mathfrak{g} = \mathfrak{h} \oplus \mathfrak{m}$. In this case, the tangent space $T_o(G/H)$ can be identified with $\mathfrak{m}$. 


By definition, the isotropy subgroup fixes the point $o$. Hence for all $h\in H$ the differential $(d\mu_h)_o:T_o(G/H)\rightarrow T_{h\cdot o}(G/H)=T_o(G/H)$ yields a linear action of $H$ on the tangent space 
called the isotropy representation $\sigma:H\longrightarrow \Gl(T_o(G/H))$, such that $\sigma (h) :=(d\mu_h)_o$. As the group acts via orientation preserving isometries, we may take the image to lie inside $\SO(T_o(G/H))$. 
We will assume furthermore that $G$ acts effectively and that therefore the isotropy representation is injective. 

By the preceding paragraphs, the adjoint representation of $G$ restricted to $H$ satisfies $\Ad^G\vert_H \cong \Ad^H \oplus \sigma$.  Therefore, if $\mathfrak{m}$ is any $H$-subrepresentation of $\mathfrak{g}$ such that $\mathfrak{g} = \mathfrak{h} \oplus \mathfrak{m}$, we have $(\Ad^G\vert_H, \mathfrak{m}) \cong (\sigma, T_o(G/H))$.


It is important to describe precisely the structure of the orthonormal frame bundle $FM$ of the homogeneous space $G/H$. By our previous considerations, this is a $\SO(\mathfrak{m})$-bundle and we can consequently identify $F_{o}M$ with $\SO(\mathfrak{m})$. Consider now the $H$-bundle $p: G\times_\sigma \SO(\mathfrak{m})\rightarrow G/H$ associated to the bundle $G\rightarrow G/H$, where the right group action of $H$ is given by $(g,A)\mapsto (gh,\sigma(h^{-1})A)=:[g,A]$. The map $$G\times_\sigma \SO(\mathfrak{m})\rightarrow FM,\quad  [g,A]\mapsto (gH,d\mu_g A)$$
is a bundle isomorphism, with $d\mu_g: F_oM\rightarrow F_{gH}M$ being the pushforward of the action by $G$.

\subsection{Classification of group actions on spheres} \label{subsection_groupactions}

We are interested in compact connected Lie groups acting effectively and transitively on spheres. Note that the compactness condition enables us to work with isometric actions, while we can always assume that the group is connected by taking the identity component, which would also act transitively on the homogeneous space. These groups were firstly classified by Montgomery and Samelson in \cite{Samuelson}. We will recall the actions of the classical Lie groups below, as we will explicitly need them later. Since a comprehensive (and comprehensible) survey of the description of these actions is somewhat hidden in the literature, we added the exceptional cases in an appendix at the end of the paper. We will also describe the loops whose class generate the fundamental group of the connected Lie groups which are not simply connected. For this purpose, we denote the usual rotation by \[R(t)=\left( \begin{matrix}
	\cos(2\pi t) & -\sin(2\pi t)\\
	\sin(2\pi t) & \cos(2\pi t)
\end{matrix}\right). \] 

\textbf{The actions of $\SO(n+1)$, $\U(n+1)$, $\SU(n+1)$ and $\Sp(n+1)$:} 

First we recall that 
 \begin{eqnarray*}
\SO(n)&:=&\{A\in \Ogroup(n):\det A=1\}, \textrm{  where }\Ogroup(n):=\{A\in \Gl(n,\R):A^tA=AA^t=Id\},\\
\SU(n)&:=&\{A\in \U(n): \det A=1\}, \textrm{ where }\U(n):=\{A\in \Gl(n,\C):A\Bar{A}^t=\Bar{A}^tA=Id\}, \\
\Sp(n)&:=&\{A\in\Gl(n,\mathbb{H}):A\Bar{A}^t=\Bar{A}^tA=Id\}.
\end{eqnarray*}
The group $\SO(n+1)$ acts transitively on the unit sphere $S^n:=\{(v_0,\dots,v_n)\in \mathbb{R}^{n+1}:\sum|v_i|^2=1\}$ by the usual matrix multiplication. It is easy to see that the isotropy group at $(1,0, \dots ,0)^T$ is given in block from by the matrix $\left(\begin{array}{c|c}
 1 & 0 \\ \hline
    0 & \SO(n)
\end{array}\right)
$. 
In order to find the generator of $\pi_1(\SO(n))=\Z_2$ for $n>2$, we use the fact that $S^n$ is $2$-connected for $n>2$ ($\pi_1(S^n)=\pi_2(S^n)=0$), together with the well-known long exact sequence of homotopy groups for the principal bundle associated to the homogeneous space structure, to conclude that the group morphism induced at the level of fundamental groups by the inclusion of the isotropy subgroup is an isomorphism (or surjective when $n=2$). Hence we can track back the generator of $\pi_1(\SO(n))$ to $\pi_1(\SO(2))$. Thus, we can write a generating loop $\alpha_n:I\longrightarrow \SO(n)$ of $\pi_1(\SO(n))$ by $\alpha_n(t)=\left(\begin{array}{c|c}	Id_{n-2} & 0 \\ \hline	0 & R(t)\end{array}\right)$.

Similarly $\U(n+1)$ and $\SU(n+1)$ act by matrix multiplication on the unit sphere in the complex $(n+1)$-dimensional space, $S^{2n+1}:=\{(v_0,...,v_n)\in \mathbb{C}^{n+1}:\sum|v_i|^2=1\}$, where $\mathbb{C}^{n+1}$ is identified with $\R^{2n+2}$, and $\Sp(n+1)$ acts on the unit sphere in the quaternions $S^{4n+3}=\{(v_0,...,v_n)\in \HH^{n+1}:\sum|v_i|^2=1\}$. Here $v_i=a_i+\ii b_i+\jj c_i+ \kk d_i$ and $|v_i|^2=a_i^2+b_i^2+c_i^2+d_i^2$ and again we can identify $\HH^{n+1}$ with $\R^{4n+4}$. The isotropy groups can be computed as above and are $\U(n)$, $\SU(n)$ and $\Sp(n)$ respectively. 

We can use the same argument as above to find a loop whose class generates $\pi_1(\U(n))=\Z$. Knowing that the loop $\beta_1:I\longrightarrow\U(1)$ such that $\beta_1(t)=e^{2\pi i t}$ generates $\pi_1(\U(1))=\Z$ it is straightforward to see that the loop we seek is $\beta_n(t)=\left(\begin{array}{c|c}	Id_{n-1} & 0 \\ \hline	0 & e^{2\pi i t}\end{array}\right)$.
    
\textbf{The actions of $\Sp( n+1)\cdot \U( 1)$ and $\Sp( n+1)\cdot \Sp( 1)$:}

The case $\Sp( n+1)\cdot \U( 1)$ and $\Sp( n+1)\cdot \Sp( 1)$ are analogous, so we will only explain the case $\Sp( n+1)\cdot \U( 1)$ in detail.

We have that $\Sp( n+1)\cdot \U( 1)=\Sp( n+1)\times \U( 1)/\{\pm(Id,1)\}$, so the group is doubly covered by $\Sp( n+1)\times \U( 1)$. We denote the elements of $\Sp( n+1)\times \U( 1)$ by $(A,z)$ and the elements of $\Sp( n+1)\cdot \U( 1)$ by $[A,z]$. Finally, we need an embedding of $\U( 1)\subset \C$ in $\Sp( 1)\subset \HH$. We choose $\iota:\U( 1)\hookrightarrow \Sp( 1)$ where $\iota(a+ib)=a+\ii b+\jj 0+ \kk 0$. Note that there are other possible inclusions. 

We define the action of $\Sp( n+1)\cdot \U( 1)$ in the following way. Given $[A,z]\in \Sp( n+1)\cdot \U( 1)$, we define $\mu_{[A,z]}:S^{4n+3}\longrightarrow S^{4n+3}$ such that $\mu_{[A,z]}(v)=Av(\iota(z))^{-1}$, where $v(\iota(z))^{-1}=(v_0(\iota(z))^{-1},...,v_n(\iota(z))^{-1})$. Note that this shows why we take $\Sp( n+1)\cdot \U( 1)$ instead of $\Sp( n+1)\times \U( 1)$ to get an effective action. The action is moreover transitive, since the usual action of $\Sp(n+1)$ on $S^{4n+3}$ is transitive. Now, we compute the isotropy subgroup $H$ for $p=(1,0,...,0)$. If $[A,z]\in H$ then $Ap=(\iota(z),0,...,0)$. From there, it is straightforward to see that  
\[H=\{[\left(\begin{array}{c|c}
   \iota(z)  & 0 \\ \hline
    0 & A
\end{array}\right),z]:A\in \Sp( n),z\in \U( 1)\}\cong \Sp( n)\cdot \U( 1).
\]
Moreover, the inclusion $f_n:\Sp( n)\cdot \U( 1)\longrightarrow \Sp( n+1)\cdot \U( 1)$ fulfills
\[f_n([A,z])=[\left(\begin{array}{c|c}
   \iota(z)  & 0 \\ \hline
    0 & A
\end{array}\right),z].
\]
Once again, all the inclusions of isotropy subgroups are isomorphisms at the level of fundamental groups. Thus we can track back the generator of $\pi_1(\Sp( n)\cdot \U( 1))=\Z$ to the case $n=0$. Then, we have that the generator the fundamental group of $\U( 1)/\Z_2$ is represented by the loop $\gamma_0:I\longrightarrow \U( 1)/\Z_2$ such that $\gamma_0(t)=[e^{i\pi t}]$. Using the inclusions, we can see that the generator of $\pi_1(\Sp( n)\cdot \U( 1))=\Z$ is represented by a loop $\gamma_n:I\longrightarrow \Sp( n)\cdot \U( 1)$ such that $\gamma_n(t)=[\iota(e^{i\pi t})Id,e^{i\pi t}]$.

The action of $\Sp( n+1)\cdot \Sp( 1)=\Sp( n+1)\times \Sp( 1)/\{\pm(Id,1)\}$ on $S^{4n+3}$, as well as its isotropy group, are the same as  $\Sp( n)\cdot \U( 1)$, just we replace the inclusion $\iota(z)$ by $z$ and take $z$ in $\Sp(1)$. A loop whose class generates the fundamental group of $\pi_1(\Sp(n)\cdot\Sp(1))\cong \Z_2 $ is $\gamma_n'(t)=[\iota(e^{i\pi t})Id,\iota(e^{i\pi t})]$.

\subsection{Spin structures on homogeneous spaces}

Let $M$ be an orientable smooth Riemannian manifold and let $(FM,p,M,\SO(n))$ be its orthonormal frame bundle. Recall that the group $\Spin(n)$ is the unique connected double covering of $\SO(n)$, which we denote by $\lambda$. Then a spin structure is a pair $(P, \Lambda)$, where $\Lambda: P \rightarrow FM$ is a 2-covering such that we have a principal $\Spin(n)$-bundle $(P,p'=p\circ\Lambda,M,\Spin(n))$ and the following diagram commutes;
\[
\begin{tikzcd}
P \times \Spin(n) \arrow[r, "\Phi_{\Spin}"] \arrow[dd, "\Lambda \times \lambda"] & P \arrow[dd, "\Lambda"] \arrow[dr, "p"] \\
 & & X \\
FM \times \SO(n) \arrow[r, "\Phi_{\SO}"] & FM \arrow [ur, "p"]
\end{tikzcd}
\]
here $\Phi_{\Spin}$ and $\Phi_{\SO}$ denote the action of these groups on the total space of their respective principal bundles. In addition, two spin structures $(P_1,\Lambda_1)$ and $(P_2,\Lambda_2)$ are said to be equivalent if there exists a principal bundle isomorphism $f:P_1\longrightarrow P_2$ such that $\Lambda_2\circ f=\Lambda_1$.    

Although we fix a Riemannian metric to define a spin structure, it is well-known that the obstruction to its existence is purely topological. More precisely, an orientable Riemannian manifold $M$ is spin if and only if its second Stiefel-Whitney class vanishes, $w_2(M)=0$. Then, there exists a one to one correspondence between spin structures up to equivalence and cohomology classes $c\in H^1(FM,\Z_2)$ such that $i^*(c)\neq 0$, where $i$ is the inclusion of a fiber. Moreover, the number of spin structures is precisely $|H^1(M,\Z_2)|$ (see \cite[Chapter 2]{Friedrich}). It is immediate to conclude from the above discussion the well known fact that $S^n$ is spin, with a unique spin structure.

The first step to study the relation between spin structures and group actions is to understand how an orientation preserving isometry $f:M\longrightarrow M$ transforms spin structures on a spin manifold $M$. Recall that each orientation preserving isometry induces an isomorphism on the frame bundle of the manifold, which we also denote by $f:FM\longrightarrow FM$, such that $f(p,(v_1,...,v_n))=(f(p),(df_{|p}v_1,...,df_{|p}v_n))$.

We would like to extend this map to a map to a spin structure $(P,\Lambda)$. More precisely, we seek a map $\tilde{f}:P\longrightarrow P$ such that 
\[
\begin{tikzcd}
	P\ar{r}{\tilde{f}}\ar{d}{\Lambda} & P\ar{d}{\Lambda}\\
	FM\ar{r}{f} & FM\\
\end{tikzcd}
\]
that must also preserve the fibers and be equivariant with respect $\Spin(n)$ right action on $P$. Since the spin structure is represented by a cohomology class $c\in H^1(FM,\Z_2)$ such that $i^*(c)\neq0\in H^1(\SO(n),\Z_2)$, then it can be seen that the map $f:FM\longrightarrow FM$ can be lifted to a map $\tilde{f}:P\longrightarrow P$ between structures if and only if $f^*(c)=c$ by using \v{C}ech cohomology (see \cite{Chichilnisky1996351}). It is important to remark that there are always two possible lifts.

Now, we can introduce the concept of $G$-invariant spin structure. Let $G$ be Lie group acting by orientation preserving isometries on a spin manifold $M$ and let $\phi_g:M\longrightarrow M$ be the isometry induced by the action of an element $g\in G$ (we use $\phi$ instead of $\mu$ to remark that the action is not necessarily transitive). As above, we also denote by $\phi_g$ the unique bundle isomorphism induced in $FM$. Then:  

\begin{defn0}\label{Ginvariant spin structure def}
We say that a spin structure $(P, \Lambda)$ is $G$-invariant if it is equipped with an action of $G$ on $P$ covering the action of $G$ on $FM$. 
\end{defn0}
This implies that for each $g \in G$ there is an isomorphism $\tilde{\phi}_g:P\longrightarrow P$ such that $\Lambda\circ\tilde{\phi}_g=\phi_g\circ\Lambda$.

\begin{remark0}\label{invariant and subgroups}
Let $G'\subset G$ be Lie groups acting on a spin manifold $M$. If a fixed spin structure is $G$-invariant, then it is also $G'$-invariant. Thus, if the spin structure is not $G'$ invariant, then it is not $G$-invariant.  
\end{remark0}

Since the action of a connected Lie group is homotopically trivial, the above cohomological condition is always fulfilled. This implies that if a connected Lie group acts on a spin manifold $M$, then either $G$ or a $2$-covering $\tilde{G}$ acts in a way that preserves the spin structure (see \cite{Chichilnisky1996351}). Note that the action of $\tilde{G}$ on $M$ is not effective. Finally, if $G$ is simply connected the only possibility for $\tilde{G}$ is $G\times \Z_2$. This implies that every spin structure is $G$-invariant when $G$ is simply connected.      

For the case of Riemannian homogeneous space, the action of $G$ on its frame bundle $G\times_\sigma\SO(n)$ is simply $\mu_g([g',A])=[gg',A]$. Assume now that there exists a lift of the isotropy representation to the spin group, that is a group morphism $\tilde{\sigma}:H\longrightarrow \Spin(n)$ such that $\sigma=\lambda\circ \tilde{\sigma}$. Note that by using covering space theory and the fact that $\Spin(n)$ is simply connected for $n>2$, it is straightforward to see that this lift exists if and only if $\sigma_*:\pi_1(H)\longrightarrow \pi_1(\SO(n))$ is trivial. Then, it is a known fact that we can construct a spin structure of the form $G\times_{\tilde{\sigma}} \Spin(n)$ with a 2-covering $\Lambda([g,x])=[g,\lambda(x)]$. The next proposition shows that the existence of a lift and the existence of a $G$-invariant spin structure are equivalent. One of the implications is well-known (see  \cite{ALEKSEEVSKYD,CGT}), while the converse statement and its proof are not clearly presented in the literature (for example, it is only mentioned in \cite{bar1992dirac}). Since it is a key step to prove the main theorem we provide a detailed proof of this fact.

\begin{prop0}\label{invariant spin structures and lifts}
Let $M=G/H$ be a spin homogeneous space where $H$ is connected. Then there exists a $G$-invariant spin structure if and only if the isotropy representation lifts to a map $\tilde{\sigma}:H\longrightarrow \Spin(n)$.
\end{prop0}

\begin{proof}
Assume that the isotropy representation lifts to a group morphism $\tilde{\sigma}:H\longrightarrow \Spin(n)$, then $\Lambda:G\times_{\tilde{\sigma}}\Spin(n)\longrightarrow G\times_\sigma \SO(n)$ is a spin structure. We consider the action of $G$ on $G\times_{\tilde{\sigma}}\Spin(n)$ such that $\tilde{\mu}_g:G\times_{\tilde{\sigma}}\Spin(n)\longrightarrow G\times_{\tilde{\sigma}}\Spin(n)$ fulfil that $\tilde{\mu}_g([g',x])=[gg',x]$ for every $g\in G$. It is clear that $\Lambda\circ\tilde{\mu}_g=\mu_g\circ\Lambda $, hence the spin structure is $G$-invariant.

Conversely suppose that we have a $G$-invariant spin structure, with total space $P \to G/H$.  Then $P$  admits a transitive action of $G \times \Spin(n)$, with $\Spin(n)$ the usual action on the fibres, and $G$ acting by the equivariant structure. Let $p \in P$ be a point and $x \in G/H$ its projection. The first projection $G \times \Spin(n) \to G$ induces an isomorphism from the isotropy $H_p$ of $p$ in $G \times \Spin(n)$ to the isotropy group $H_x$. This is because $\Spin(n)$ acts simply transitively on the fibres of $P \to G/H$, so that  for every element $h \in H_x$ fixing $x$, there must be a unique lift $h \in H_p$ to an element fixing $p$.  The inverse isomorphism $H_x \to H_p \subseteq G \times \Spin(n)$ gives, by projection to the second factor, a homomorphism $H_x \to \Spin(n)$. It is clear that, when we project this further to $\SO(n)$, we obtain the isotropy action on the frame bundle at $x$. \qedhere

\end{proof}

If the lift exists then it is unique since $H$ is connected. In this case, a $G$-invariant spin  structure  $(P,\Lambda)$ on $G/H$ is equivalent to the spin structure $G\times_{\tilde{\sigma}}\Spin(n)$ given by the map $f:G\times_{\tilde{\sigma}}\Spin(n)\longrightarrow P$ such that $f([g,x])=\tilde{\mu}_g(p_ox)$, where $\Lambda(p_o)=[e,Id]$. This leads to the next corollary: 
\begin{corollary0}
If there exists a $G$-invariant spin structure on $G/H$ with $H$ connected, then it is unique.
\end{corollary0}

\begin{remark0}
Let $G$ be a connected Lie group, seen as the homogeneous spaces  $G\cong G/\{e\}$. Then, the isotropy representation lifts trivially and we have a $G$-invariant spin structure, which is corresponds to the trivial bundle $G\times \Spin(n)$. This spin structure is the only $G$-invariant spin structure on $G$. This fact is especially relevant when $G$ has multiple inequivalent spin structures, like $T^n$, since it gives a preferred spin structure to work with. 
\end{remark0}

In the spirit of the main theorem, we can also ask which connected Lie groups $G$ acting transitively on other Lie groups $M$, seen as a manifold, leave its trivial spin structure invariant. There are cases, like $M=T^n$, where the only group acting transitively on it is the group itself, thus the answer is trivial. However, if we pick $M=\SU(2)\cong S^3$, then $\SO(4)$ also acts transitively and effectively on $M$, but as we will see later, the unique spin structure of $S^3$ is not $\SO(4)$-invariant. 

Now as we mentioned above, one of the interesting properties of the Lie groups acting transitively and effectively on spheres, is that they are almost in one to one correspondence with the Riemannian holonomy groups. In fact we have the following nice property, which seems as well as its proof to be missing from literature.

\begin{proposition}\label{equivalence holonomy}
Let $G$ be the holonomy group of a simply-connected irreducible non-symmetric Riemannian manifold
of dimension $n+1 \geq 3$. Let $H<G$ be a subgroup such that $G/H \cong S^n$ given by Berger's classification.
Then there exists a lift of the holonomy representation $h:G\rightarrow \SO(n+1)$ to $\Spin(n+1)$ if and and only if $S^n$ has a $G$-invariant spin structure.
\end{proposition}
\begin{proof}
By Proposition \ref{invariant spin structures and lifts} the datum of a $G$-invariant spin structure of the sphere $G/H$ is equivalent to a lift of the isotropy representation $\tilde{\sigma}:H\longrightarrow \Spin(n)$. Now we consider the following commutative diagram of fundamental groups 
\[
\begin{tikzcd}
\arrow[dotted]{d}& \arrow[dotted]{d}\\
	\pi_1(H)\ar{r}{\sigma_*}\ar{d}{} & \pi_1(\SO(n))\ar{d}{}\\
	\pi_1(G)\ar{r}{h_*} \arrow[dotted]{d}& \pi_1(\SO(n+1))\arrow[dotted]{d}\\
 ~&~\\
\end{tikzcd}
\]
where the vertical arrows come from the long exact sequences of homotopy groups for the quotient space $G/H$ (resp. $\SO(n+1)/SO(n)$). Since these spaces are spheres of dimension greater than two, the two vertical arrows are isomorphisms, which finishes the proof, since then the map $\sigma_*$ is trivial if and only if $h_*$ is trivial by above mentioned criterion for the existence of a lift. 
\end{proof}
\section{Main Result}



The groups $\SU(n)$, $\Sp(n)$, $\eG$, $\Spin(7)$ and $\Spin(9)$ are simply connected, so in these cases the isotropy representations always lift and the spin structure is invariant. Hence, we only need to consider the cases $\SO(n)$, $\U(n)$, $\Sp(n)\cdot\U(1)$ and $\Sp(n)\cdot\Sp(1)$. 

We present two different approaches to prove the theorem for these groups. Our first approach is to compute $\sigma(\gamma(t))=(d\mu_{\gamma(t)})_o$, where $\gamma(t)$ is a loop whose class generates the fundamental group of the group we study. 
In the second approach we use representation theory to find $\sigma$ by using that $\Ad^G\vert_H =\Ad^H \oplus \sigma$, and then compute $\sigma(\gamma(t))$. 

For each of the following lemmata corresponding to the four cases above, we provide first the proof using the first approach and then the more representation theoretical technique. Note that for the complex and quaternionic group, we will be using complex representations which are more convenient to deal with.  One can uniquely recover, up to isomorphism, the isotropy representation from its complexification: if $\rho, \rho'$ are real finite-dimensional representations of any group, then $\rho \cong \rho'$ if and only if the complexifications $(\rho \otimes \mathbb{C}) \cong (\rho' \otimes \mathbb{C})$ are isomorphic. This follows because
$(\rho \otimes \mathbb{C})^{\mathbb{R}} \cong \rho \otimes (\mathbb{C})^{\mathbb{R}} \cong \rho \oplus \rho$. Here $\psi^{\mathbb{R}}$ denotes the underlying real representation of a complex one $\psi$.

\begin{lemma0}\label{SO(n) not invariant}
The spin structure of $S^n$ is not $\SO(n+1)$-invariant.
\end{lemma0}

\begin{proof}
Recall that in this case, $S^n=\SO(n+1)/\SO(n)$ and we take $o=(1,0,...,0)$. Then the action is by restriction of a linear action, hence the differential $(d\mu_{\alpha_n(t)})_o:T_o\R^{n+1}\longrightarrow T_o\R^{n+1}$ is the matrix $f_n(\alpha_n(t))$, where $\alpha_n$ is the generating loop described in Section \ref{subsection_groupactions}. If we restrict it to  $T_pS^{n}=\{v\in\R^{n+1}:\langle v,(1,0,...,0) \rangle=0\}=\langle e_2,...,e_{n+1}\rangle$, we obtain that $\sigma(\alpha_n(t))=\alpha_n(t)\in \SO(n)$. Since $[\alpha_n(t)]$ is the generator of $\pi_1(\SO(n))$, we can conclude that the isotropy representation does not lift, which implies that the spin structure is not $\SO(n+1)$-invariant. 

Let us now denote by $\lambda_n$ the standard representation of $\SO(n)$, which is given by the left multiplication on $\R^n$. It is a well-known fact that $\Ad^{\SO(n)} \cong \bigwedge^2\lambda_n$. Recall that $\sigma=\Ad^{\SO(n+1)}|_{\SO(n)}:\SO(n)\longrightarrow \SO(\mathfrak{m})$ and that we have the chain of isomorphisms $\Ad^{\SO(n+1)}|_{\SO(n)}\cong \bigwedge^2(\lambda_{n+1}|_{\SO(n)})\cong \bigwedge^2(\lambda_n\oplus 1)\cong \bigwedge^2\lambda_{n}\oplus\lambda_n$, where $1$ denotes the trivial representation. The first summand is the adjoint representation of the isotropy subgroup $\SO(n)$, which means that $\sigma\cong\lambda_n$. This implies, as we have already seen, that $\sigma(\alpha_n(t))=\alpha_n(t)$ and hence the isotropy representation does not lift. 
\end{proof}

The next case we discuss is the action of $\U(n)$.

\begin{lemma0}\label{U(n) not invariant}
The spin structure of $S^{2n+1}$ is not $\U(n+1)$-invariant.
\end{lemma0}   
\begin{proof}
We proceed in the same way for $S^{2n+1}=\U(n+1)/\U(n)$. Let again $\beta_n$ be the generating loop described in Section \ref{subsection_groupactions}. In order to compute $\sigma(\beta_n(t)) $, we view the action as a linear action on $S^{2n+1}\subset \C^{n+1}\subset\R^{2n+2}$. We can decomplexify the matrices of $\U(n+1)$ in order to see them in $\SO(2n+2)$. Consequently, the isotropy representation is given by the natural inclusions $ \U(n)\subset \SO(2n)\subset \SO(2n+1)$. Therefore, we obtain that $\sigma(\beta_n(t))=\left(\begin{array}{c|c}	Id_{2n-1} & 0 \\ \hline	0 & R(t)\end{array}\right)\in \SO(2n+1)$. Like in the above case, this means that the isotropy representation does not lift and the spin structure is not $\U(n)$-invariant.

We now look at the second approach. Let $\mu_n$ be the standard complex representation of $\U(n)$ as a matrix  acting on $\C^{n}$. Then, the complexified adjoint representation satisfies $\Ad^{\U(n+1)}\otimes \C \cong \mu_{n+1} \otimes_{\C} \mu_{n+1}^*\cong \mu_{n+1} \otimes_{\C} \bar{\mu}_{n+1}$. If we restrict it to $\U(n)$, we obtain that $\Ad^{\U(n+1)}\otimes \C \vert_{\U(n)} \cong (\mu_n \oplus 1) \otimes_{\C} (\bar{\mu}_n \oplus 1)	\cong (\mu_n \otimes \bar{\mu}_n) \oplus \bar{\mu}_n \oplus \mu_n \oplus 1 $. Therefore, the complexified isotropy representation is isomorphic to $\bar{\mu}_n \oplus \mu_n \oplus 1$.  Note that this is isomorphic to the complexification of $\mu_n^{\mathbb{R}} \oplus 1$ (where now $1$ is the real trivial representation). So the real isotropy representation is isomorphic to this.

If we apply the isotropy representation to $\beta_n(t)$, we obtain $[R(t) \oplus 1 \oplus...\oplus 1 ]$.
So we will obtain a rotation matrix that generates the fundamental group of $\SO(2n+1)$.
\end{proof}

For the sake of completeness, we also compute the (complexified) isotropy representation for the Lie groups $\SU(n)$ and $\Sp(n)$. 
In the first case, the sphere is the same as for $\U(n)$, with a restricted action, so we simply restrict the previous isotropy representation to $\SU(n)$.

For the second case, let $\nu_n$ be the standard complex representation of $\Sp(n)$ (of complex dimension $2n$). Then, $\Ad^{\Sp(n+1)}\otimes \C=S^2(\nu_{n+1})$, where $S^2$ denotes the second complex linear symmetric power. Thus, $\Ad^{\Sp(n+1)}\otimes \C|_{\Sp(n)}\cong S^2(\nu_{n}\oplus1 \oplus 1)\cong S^2(\nu_n)\oplus(S^1(\nu_n)\otimes S^1(1\oplus 1))\oplus S^2(1 \oplus 1)\cong S^2(\nu_n)\oplus\nu_n\oplus \nu_n \oplus 1 \oplus 1 \oplus 1$. This implies that the complexified isotropy representation is isomorphic to the last sum after removing the first summand. Note that this is isomorphic to the complexification of $\nu_n^{\mathbb{R}} \oplus 1 \oplus 1 \oplus 1$. Thus the real isotropy representation is isomorphic to this.

\begin{lemma0}\label{Sp(n)U(1) Sp(n)Sp(1) invariant}
The spin structure of $S^{4n+3}$ is invariant by the transitive actions of $\Sp(n+1)\cdot\U(1)$ and $\Sp(n+1)\cdot\Sp(1)$ if and only if $n$ is odd.
\end{lemma0}  

\begin{proof}
Both case are analogous, so we focus mainly in the case $S^{4n+3}=\Sp(n+1)\cdot\U(1)/\Sp(n+1)\cdot\U(1)$. Like in the previous cases, we start by computing $\sigma([A,z])=(d\mu_{f_n[A,z]})_o$. In order to do this, we will see that $\mu_{f_n[A,z]}$ is a linear map acting on $S^{4n+3}\subset \R^{4n+4}$. We have that $\mu_{\gamma_n(t)}(v)=\iota(e^{i\pi t})v\iota(e^{i\pi t})^{-1}=(C_{\iota(e^{i\pi t})}v_0,...,C_{\iota(e^{i\pi t})}v_n)$, where $C_{w}:\Sp( 1)\longrightarrow \Sp( 1)$ is the conjugation by an element $w\in \Sp( 1)$. Then, a tedious computation shows that the conjugation induces a linear map on $\R^4$ (which we can restrict to $S^3$) given by a $4\times 4$ matrix
\[\left(\begin{array}{c|c}
	Id & 0 \\ \hline
	0 & R(t)
\end{array}\right),
\]

Therefore, we have that $\mu_{f_n(\gamma_n(t))}\in \SO(4n+4)$ is 
\[\left(\begin{array}{c|c|c}
	\begin{array}{c|c}   Id & 0 \\ \hline 0 & R(t)\end{array} & \cdots &0  \\ \hline
	\vdots & \ddots & \vdots \\ \hline
	0 & \cdots & \begin{array}{c|c}   Id & 0 \\ \hline 0 & R(t)\end{array}\\
\end{array} \right).
\]

Now, we use that $T_pS^{4n+3}=\{v\in\R^{4n+4}:\langle v,(1,0,...,0) \rangle=0\}=\langle e_2,...,e_{4n+4}\rangle$. This means that we have that 
\[
(d\mu_{f_n(\gamma_n(t))})_o=\left(\begin{array}{c|c|c}
	\begin{array}{c|c}   1 & 0 \\ \hline 0 & R(t)\end{array} & \cdots &0  \\ \hline
	\vdots & \ddots & \vdots \\ \hline
	0 & \cdots & \begin{array}{c|c}   Id & 0 \\ \hline 0 & R(t)\end{array}\\
\end{array} \right)\in \SO(4n+3).
\]
Since, we have $n+1$ rotations, we can conclude that $\sigma_*[\gamma_n]=n+1\mod 2$. Hence, the isotropy representation lifts when $n$ is odd. In other words, the spin structure of $S^{8m-1}$ is invariant by the action of $\Sp(2m+2)\cdot \U(1)$ and $\Sp(2m+2)\cdot \Sp(1)$.

In order to use our second method, we need again to find the isotropy representation first.  We note that $G=\Sp(n+1)\cdot\Sp(1)$ and $\Sp(n+1) \times \Sp(1)$ both have the same (complexified) Lie algebra $\mathfrak{g}^{\mathbb{C}}\cong\mathfrak{sp_{\mathbb{C}}}(2n+2)\oplus \mathfrak{sp_{\mathbb{C}}}(2)$ (note that here the dimension is not the quaternionic but the complex dimension). Moreover $\mathfrak{g}^{\mathbb{C}}=\mathfrak{m}^{\mathbb{C}}\oplus\mathfrak{h}^{\mathbb{C}}$, where $\mathfrak{m}^{\mathbb{C}}$ is isomorphic, as an  $\mathfrak{h}^{\mathbb{C}}$-representation, to the complexified isotropy representation.

Since by Section \ref{subsection_groupactions}, $H=\{[\left(\begin{array}{c|c}
    z & 0 \\ \hline
    0 & A
\end{array}\right),z]:A\in \Sp( n),z\in \Sp( 1)\}\cong \Sp( n)\cdot \Sp( 1)$,
 we have that \[\mathfrak{h}^{\mathbb{C}}=\{(\left(\begin{array}{c|c}
  \xi  & 0 \\ \hline
    0 & \mathfrak{a}
\end{array}\right),\xi):\mathfrak{a}\in \mathfrak{sp_{\mathbb{C}}}(2n), \xi\in \mathfrak{sp}_\mathbb{C}( 2)\}\cong \mathfrak{sp_{\mathbb{C}}}(2n)\oplus \mathfrak{sp_{\mathbb{C}}}(2)'\subseteq  \mathfrak{g}_{\mathbb{C}},\]
where $\mathfrak{sp_{\mathbb{C}}}(2)'=\{(\left(\begin{array}{c|c}
  \xi  & 0 \\ \hline
    0 & 0
\end{array}\right),\xi):\xi\in \mathfrak{sp}_\mathbb{C}( 2)\}$.

Note that for every direct sum of Lie algebras $\mathfrak{g} = \mathfrak{g}_1 \oplus \mathfrak{g}_2$, the first summand is a subrepresentation of $\mathfrak{g}$ under the adjoint action; on $\mathfrak{g}_1$, $\mathfrak{g}_1$ acts by its adjoint representation, and $\mathfrak{g}_2$ acts trivially.  
Consequently if $\mathfrak{h}$ is a Lie subalgebra of $\mathfrak{g}$, its action on $\mathfrak{g}_1$ is via the composition $\displaystyle \mathfrak{h} \to \mathfrak{g} \twoheadrightarrow \mathfrak{g}_1 \mathop{\to}^{\Ad} \Gl(\mathfrak{g}_1)$.


Under the projection to the first summand 
$\mathfrak{sp}_{\mathbb{C}}(2n+2)$, $\mathfrak{h}^{\mathbb{C}}$ maps isomorphically to $\mathfrak{sp}_{\mathbb{C}}(2) \oplus \mathfrak{sp}_{\mathbb{C}}(2n)$. By the preceding paragraph, the action of $\mathfrak{h}^{\mathbb{C}}$ on the first summand $\mathfrak{sp_{\mathbb{C}}}(2n+2)$ of $\mathfrak{g}^{\mathbb{C}}$ is the restriction of the adjoint representation of $\mathfrak{sp_{\mathbb{C}}}(2n+2)$ to $\mathfrak{sp_{\mathbb{C}}}(2n)\oplus \mathfrak{sp_{\mathbb{C}}}(2)\subseteq \mathfrak{sp_{\mathbb{C}}}(2n+2)$.


  
Now, a complement to $\mathfrak{h}^{\mathbb{C}}$ in $\mathfrak{g}^{\mathbb{C}}$ can be obtained as the complement to $\mathfrak{sp}_{\mathbb{C}}(2n)$ in the first factor $\mathfrak{sp}_{\mathbb{C}}(2n+2)$. Namely, as we computed previously,
$$\mathfrak{sp}_{\mathbb{C}}(2n+2) \cong S^2(\mathbb{C}^{2n+2}) \cong
S^2 \mathbb{C}^{2n} \oplus (\mathbb{C}^{2n} \otimes_{\mathbb{C}} \mathbb{C}^2) \oplus S^2 \mathbb{C}^2,$$
with the first summand $S^2 \mathbb{C}^{2n}$ corresponding to $\mathfrak{sp}_{\mathbb{C}}(2n)$. So we can identify $\mathfrak{m}^{\mathbb{C}}$ with $(\mathbb{C}^{2n} \otimes_{\mathbb{C}} \mathbb{C}^2) \oplus 
\mathfrak{sp}_{\mathbb{C}}(2)$. 
 In terms of matrices this is:
$$\mathfrak{m}^{\mathbb{C}} = \{(\left(\begin{array}{c|c}
   *  & * \\ \hline
    * & 0
\end{array}\right),0)\} \in \mathfrak{sp}_{\mathbb{C}}(2n+2) \oplus \mathfrak{sp}_{\mathbb{C}}(2).$$
As a result we can realize the complexified isotropy representation as $(\nu_n \boxtimes_{\mathbb{C}} \nu_1) \oplus (\Ad^{\Sp(1)} \otimes \mathbb{C})$, where $\nu_n$ is the standard representation of $\Sp(n)$.  Here the action on the second summand is via the projection $H=\Sp(n) \cdot \Sp(1) \to \Sp(1)/\pm I$.  This representation is isomorphic to the complexification of $\tilde{\nu}_n^{\mathbb{R}} \oplus \Ad^{\Sp(1)}$. Here, $\tilde{\nu}_n^{\mathbb{R}}$ denotes the real representation corresponding to the first summand, which is nothing but $\mathbb{H}^n$ with the action described in Section \ref{subsection_groupactions} (note that $\Sp(1)$ does not act trivially). Note that for the second summand, the action of $\Sp(n)\cdot\Sp(1)$ factors through the projection $\Sp(n)\cdot\Sp(1) \twoheadrightarrow \Sp(1)/\pm I$ (whose adjoint representation is the same as that of $\Sp(1)$ itself).



In the case of $\Sp(n)\cdot\U(1)$, including $\U(1)$ into $\Sp(1)$ as above, we must restrict the factor for $\Sp(1)$ to that of $\U(1)$. This yields the complexified representation $\nu_n \boxtimes \nu_1 \vert_{\U(1)} \oplus \Ad^{\Sp(1)}\vert_{\U(1)}\otimes \mathbb{C})$.

Note that $\nu_1|_{\U(1)}$ and $\Ad^{\Sp(1)}|_{\U(1)} \otimes \mathbb{C}$ are reducible, as they are complex representations of the abelian group $\U(1)$ of finite dimension greater than one. It is well-known how to decompose these into one-dimensional representations (e.g., by viewing $\U(1)$ as the maximal torus inside $\Sp(1)$ and then decomposing the standard and  complex adjoint representations of $\Sp(1)$, or alternatively of the complexified Lie algebra, isomorphic to $\mathfrak{sl}_{\mathbb{C}}(2)$).
The one-dimensional representations of $\U(1)$ are of the form $\rho_m: z \mapsto (z^m)$ for $m \in \mathbb{Z}$.  We have $\nu_1|_{\U(1)} \cong \rho_1 \oplus \rho_{-1}$ and $\Ad^{\Sp(1)}|_{\U(1)} \otimes \mathbb{C} \cong \rho_{2} \oplus \rho_0 \oplus \rho_{-2}$, where $\rho_0 = 1$ is the trivial representation. Putting this together, the complexified isotropy representation decomposes as:


\[\nu_n \boxtimes (\rho_{-1}\oplus\rho_1) \oplus\rho_2 \oplus 1 \oplus \rho_{-2}.\]

This is isomorphic to the complexification of $\tilde{\nu}_n^{\mathbb{R}}|_{\Sp(n) \cdot U(1)} \oplus \rho_2^{\mathbb{R}} \oplus 1$. Here $\rho_2^{\mathbb{R}}$ is a representation of $U(1)$ which 
factors through the quotient $\U(1)/\pm I$, so its overall representation of $\Sp(n) \cdot U(1)$ is given as the composition $\displaystyle \Sp(n) \cdot \U(1) \twoheadrightarrow \U(1)/\pm I \mathop{\to}^{\rho_2} \Gl(1,\mathbb{C})$.
So this sum is isomorphic to the real isotropy representation.

Plugging in the generator calculated above, one obtains:

\[\left(\begin{array}{c|c|c}
	\begin{array}{c|c}   Id & 0 \\ \hline 0 & R(t)\end{array} & \cdots &0  \\ \hline
	\vdots & \ddots & \vdots \\ \hline
	0 & \cdots & \begin{array}{c|c}   Id & 0 \\ \hline 0 & R(t)\end{array}\\
\end{array} \right)
\oplus (R(t)) \oplus 1.
\]




%
%
We note that this matrix contains $n+1$ generators of the fundamental group of $\SO(4n+3)$ on the diagonal, and can thus be expressed as the product of $n+1$ generators of the fundamental group $\mathbb{Z}/2$. This is equal to the parity of $n+1$, so the isotropy representation lifts for $n+1$ even.
\end{proof}


We summarise the information of the transitive actions on spheres in Table \ref{table:isotropy}.
\begin{table}[h!] 
\centering
\begin{tabular}{ |c||c|c|c|c|  }
	\hline
	Lie group & Manifold & Isotropy Subgroup &Isotropy Representation & $G$-invariant spin ?\\
	\hline
	$\SO(n+1)$   &$S^{n}$    & $\SO(n)$&   $\lambda_n$ & No\\
	$\U(n+1)$   &$S^{2n+1}$    & $\U(n)$&   $\mu_n^{\mathbb{R}} \oplus 1$ & No\\
	$\SU(n+1)$   &$S^{2n+1}$    & $\SU(n)$&   $\mu_n^{\mathbb{R}} \oplus 1$ & Yes\\
	$\Sp(n+1)$   &$S^{4n+3}$    & $\Sp(n)$&   $\nu_n^{\mathbb{R}} \oplus 1 \oplus 1 \oplus 1$ & Yes\\
	$\Sp(n+1)\Sp(1)$   &$S^{4n+3}$    & $\Sp(n)\Sp(1)$& $\tilde{\nu}_{n}^{\mathbb{R}} \oplus \Ad^{\Sp(1)}$ & $n$ odd\\
	$\Sp(n+1)\U(1)$   &$S^{4n+3}$    & $\Sp(n)\U(1)$ &  $\tilde{\nu}_{n}^{\mathbb{R}} 
\oplus 
\Ad^{\Sp(1)}\vert_{\U(1)} $ 
& $n$ odd\\
	$\eG$   &$S^6$    & $\SU(3)$&   $\mu_3^{\mathbb{R}}$ & Yes\\
	$\Spin(7)$   &$S^7$    & $\eG$&  $\zeta$ & Yes\\
	$\Spin(9)$   &$S^{15}$    & $\Spin(7)$&  $\lambda_7\oplus \Delta_7$ & Yes\\
	\hline
\end{tabular}
\begin{center}
\caption{Invariance of the spin structure on the sphere by transitive and effective compact Lie group actions.}\label{table:isotropy}	
\end{center}
\end{table}

\subsection{Closing remarks}

In this section we discuss several applications of our results.

A first immediate consequence of the main theorem is a very short and basic proof of the  following fact proved in \cite[X.\S 21]{Marchifava_1975} and \cite[Proposition 2.3]{Salamon1982}. 
\begin{corollary}
	Any $8k$-dimensional quaternionic K\"ahler manifold is spin.
\end{corollary}
\begin{proof}
	It is well-known that the holonomy group of a quaternionic K\"ahler manifold $M$ of real dimension $4n+4$ is $\Sp(n+1)\cdot \Sp(1)$. But then by Proposition \ref{equivalence holonomy}, we know that there exists a lift of the holonomy representation to the spin group if and only if there exists a lift of the isotropy representation $H=\Sp(n)\cdot \Sp(1)$ for the corresponding sphere $S^{4n+3}$. This happens exactly in the case where $n$ is odd by lemma \ref{Sp(n)U(1) Sp(n)Sp(1) invariant}, which proves our claim.
\end{proof}
Note that by the same type of argument we also get immediately the well-known facts that Calabi-Yau manifolds (holonomy $\SU(n)$), hyperk\"ahler manifolds (holonomy $\Sp(n)$) and $\eG$ and $\Spin(7)$ -manifolds are spin, though these facts are obvious since in these cases $G$ is simply connected. We want to point out that although the group $\Sp(n)\cdot \U(1)$ is not a Riemannian holonomy group, it is nevertheless linked to some interesting geometric structures as it is appearing as a weak holonomy group for particular classes of connections with torsion (see \cite{Alexandrov2006}). 

Now, as we have shown in the main theorem, there are cases where the spin structure of the sphere is not $G$-invariant, but we can then take a double covering which does preserve the spin structure. In all of the non-lifting cases, there is in fact a unique connected double covering up to isomorphism.  The next remark focuses on them.

\begin{remark0}
For $\SO(n+1)$ acting on $S^n$, the double covering is precisely $\Spin(n)$ acting non-effectively  on $S^n$. More explicitly, the action is given by $\mu_{x}(v)=\lambda(x)v$ for $x\in \Spin(n+1)$ and $v\in S^{n}$. 	

For $\U(n+1)$ acting on $S^{2n+1}$, the double covering is $\operatorname{MU}(n+1) := \{(A, z) \in U(n+1) \times \mathbb{C}^\times \mid \det A = z^2\}$ (sometimes called the metaunitary group).

For $\Sp(n+1)\cdot \Sp( 1)$ and $\Sp(n+1) \cdot \U(1)$ the double coverings are $\Sp(n+1) \times \Sp(1)$ and $\Sp(n+1) \times \U(1)$, respectively.


\end{remark0}
Now, given any connected group $G$ acting transitively on a sphere by isometries via $\alpha: G \to \operatorname{Isom}(S^n)$, then $G/\ker(\alpha)$ is one of the groups in Table \ref{table:isotropy}. Then there exists a $G$-invariant spin structure on $S^n$ if and only if either the last column is ``Yes'' or the map $G \twoheadrightarrow G/\ker(\alpha)$ factors through the aforementioned connected double covering. This classifies all groups $G$ acting transitively on a sphere by isometries for which the sphere admits a $G$-invariant spin structure.

The next proposition describes a nice relationship between $G$-invariant spin structures on $G/H$ and finite subgroups of $G$:

\begin{prop0}
Let $M=G/H$ be a spin homogeneous space where $H$ is connected. Then, a spin structure in $G/H$ is $G$-invariant if and only if the spin structure is invariant by all subgroups of $G$ isomorphic to $\Z_2$.
\end{prop0}

\begin{proof}
One direction is trivial, since if the spin structure is $G$-invariant, then it is also invariant by all its subgroups. Thus, we only need to proof the converse. Assume that there is a spin structure which is not $G$-invariant, then there exists $\alpha\in\pi_1(H)$ such that $\sigma_*(\alpha)\neq 0$. Firstly, we want to find an adequate loop $a:S^1\longrightarrow H$, such that $[a]=\alpha$. Let $i:T^r\hookrightarrow H$ be the inclusion of a maximal torus, then the induced map $i_*:\Z^r\cong \pi_1(T,e)\longrightarrow\pi_1(H,e)$ is surjective \cite[Theorem 7.1]{brocker2013representations}. Alternatively, one can see that $G/T$ is simply connected by arguing that $G/T$ is a complex manifold with a CW-structure with only even dimensional cells, and then use the long exact sequence of homotopy.     
	
Then, there exists $x=(x_1,...,x_r)\in \Z^r$ such that $i_*(x)=\alpha$. Note that we may assume without lost of generality that $x$ is primitive in $\Z^r$ (i.e there is no $y\in\Z^r$ and $\lambda\in\Z$ with $|\lambda|>1$ such that $x=\lambda y$). Thus, we have a map $\tilde{a}:\R\longrightarrow \R^r$ such that $\tilde{a}(t)=xt$, which induces a loop $a:S^1\longrightarrow T^s$ whose class is $(x_1,...,x_r)$. Moreover, this map is an injective group morphism. Therefore, we have seen that we can choose the representative of $\alpha\in\pi_1(H)$ to be an injective group morphism $a:S^1\longrightarrow H$ (note that if we choose a non primitive element, the map may not be injective). 
	
We want to see that the involution $a(\tfrac{1}{2})$ does not preserve the spin structure. Because $\sigma_*(\alpha)\neq 0$, we can deduce that the spin structure is not $S^1$-invariant, thus we have the short exact sequence
$$ 1\longrightarrow \Z_2\longrightarrow S'\longrightarrow S^1\longrightarrow 1,$$ where $S'\cong S^1$ is the double covering which preserves the spin structure of $G/H$. Then, the lift of $\Z_2$ is a subgroup of $S^1$, which implies that it is $\Z_4$. Consequently, the involution action $\Z_2=\{Id,a(\tfrac{1}{2})\}$ does not leave the spin structure invariant.
\end{proof}

\begin{example0}
We can find explicitly a subgroup $\Z_2$ in the case of $S^n$ with the group action of $\SO(n)$. By the above proposition, we can pick the involution $f=\mu_{\alpha_n(\tfrac{1}{2})}$. The isometry $f:S^n\longrightarrow S^n$ is given by restricting the linear map $A:\R^{n+1}\longrightarrow \R^{n+1}$, where $A\in \SO(n+1)$ is the matrix
\[A=\left(\begin{array}{c|c}	Id_{n-1} & 0 \\ \hline	0 & -Id_2\end{array}\right).\]

We claim that the spin structure of $S^n$ is not invariant under this action. Indeed, the induced map by $f$ on the frame bundle $FS^n=\SO(n+1)$, $f_A:\SO(n+1)\longrightarrow \SO(n+1)$, is simply the matrix multiplication $X\mapsto AX$. Thus it lifts to a map $\tilde{f}:\Spin(n+1)\longrightarrow \Spin(n+1)$ which fulfils that $\tilde{f}(x)=(e_{n} \cdot e_{n+1}) \cdot x$, where $\cdot$ is the Clifford multiplication (see the appendix for details on Clifford algebras). Since $e_n \cdot e_{n+1} \cdot e_n \cdot e_{n+1}=-1$, the map $\tilde{f}^2(x)=-x$. Hence, the lift of the action on the spin structure is given by the group $\Z_4$, which implies that the spin structure is not invariant.
\end{example0}

Finally, it would be interesting to study the same question of G-invariant spin structures on other homogeneous spaces. The starting point of this problem has two difficult conditions to check; firstly,given a homogeneous space $G/H$, we need to know which Lie groups acts transitively and effectively on it and secondly  we need to know if it admits a spin structure. The first condition has been studied in \cite{HsiangSu1968,Oniscik_1968}, while we refer to \cite{CGT,ALEKSEEVSKYD} for the study of spin structures on some classes of homogeneous spaces. 

For example, we can consider the complex projective space $\CP^{m}$, which is known to be spin if and only if $m$ is odd. Then, if $m=2n+1$, the only groups that act effectively and transitively on them are $\SU(2n+2)$ and $\Sp(n+1)$ (see \cite{Oniscik_1968}). Therefore, its unique spin structure is invariant by both group actions.	

\appendix
\section*{Appendix: Exceptional group actions and isotropy representations}
The purpose of this appendix is to give a self-contained, simple and concise description of the three exceptional Lie group actions on spheres based on the octonions, and especially the realization of $G_2$ as an isotropy group, which have been scattered or missing in the literature so far. For some of the facts mentioned here for $S^6$ and $S^7$ we refer to \cite{Draper2018}.

Since the octonions are a crucial tool in understanding the three exceptional Lie group actions, we will recall their properties quickly. Let $\OO$ be the eight dimensional octonion algebra over $\mathbb{R}$, with basis $\{1=i_0, i_1,\dots, i_7\}$  satisfying the usual relations:
\begin{eqnarray*}
1i_k&=&i_k1=i_k, \textrm{ for all k in }\{0,\dots 7\}\\
i_k^2&=&-1, \textrm{ for all k in} \{1,\dots 7\}\\
i_ki_{k+1}&=&i_{k+3},\, i_{k+1}i_{k+3}=i_k,\, i_{k+3}i_k=i_{k+1},\textrm{ for all k in} \{1,\dots 7\},
\end{eqnarray*}
where the last indices are interpreted modulo 7. Any element of $\OO$ can be written as $z=\sum_{k=0}^7a_ki_k$, for some real numbers $a_k$. The octonionic conjugation is given by $\bar{z}:=a_0i_0-\sum_{k=1}^7a_ki_k$ and the real and imaginary part of $z$ by $\textrm{Re} z:=\frac{z+\bar{z}}{2}$ and $\textrm{Im} z:=\frac{z-\bar{z}}{2}$ respectively. The inner product is defined for all $x,y$ in $\OO$ by $\langle x,y\rangle=\textrm{Re}(x\bar y)$, which is the usual scalar product on $\R^8$. Any triple of unit imaginary octonions $(e_1,e_2,e_3)$, such that $e_1$ anticommutes with $e_2$, and $e_3$ anticommutes with $e_1,\, e_2$ and $e_1e_2$, known as a Cayley triple, generates the entire octonions.\\ 
\paragraph{\textbf{The action of $\eG$ on $S^6$:}} The group $\eG$  is the automorphism group $$\eG:=\{\phi\in \Gl_{\R}(\OO):\phi(xy)=\phi(x)\phi(y),\,\forall x,y\in\OO\}$$
of the 8-dimensional real vector space $\OO\cong\R^8$. Since $\phi(1)=1$, and $\phi$ preserves the inner product $\langle\cdot,\cdot\rangle$ for all $\phi$ in $\eG$, $\eG$ acts on the imaginary octonions $\textrm{Im}(\OO)\cong\R^7$, and $\eG\leq \SO(\textrm{Im}(\OO))\cong\SO(7)$. We call this seven dimensional irreducible representation $\zeta: \eG\rightarrow \SO(\textrm{Im}(\OO))$. 

The action restricts to a left action $\eG\times S^6\rightarrow S^6,\,(\phi,v)\mapsto\phi(v)$ on the unit imaginary octonions $S^6:=\{(v_0,...,v_n)\in \textrm{Im}(\OO)\cong\R^7:\sum|v_i|^2=1\}$. Since, given any two Cayley triples, there exists a unique automorphism of $\OO$ mapping the first to the second, this last action is transitive. 

We now want to compute the isotropy subgroup ${\eG}_{i_1}$ stabilizing the element $i_1$. Note that the field of complex numbers lives naturally  inside of $\OO$ as $\mathbb{C}=span\{1,i_1\}$.  Any element $z$ in $\OO$ can be written as $z=(a_0+a_1i_1)+(a_2+a_4i_1)i_2+(a_3+a_7i_1)i_3+(a_5+a_6i_1)i_5=:z_0+z_1i_2+z_2i_3+z_3i_5$, and we get an isomorphism 
$f:\OO\rightarrow\C\oplus\C^3,\,z\mapsto z_0+(z_1,z_2,z_3)^T$. It is easy to see that if $\phi\in G_2$ is an automorphism fixing $i_1$, then $\phi(z)=z_0+z_1\phi(i_2)+z_2\phi(i_3)+z_3\phi(i_5)$. It therefore fixes $\C$, preserves $\C^3$, and induces a $\C$-linear transformation of $\C^3$ given by $f\circ\phi(z)=z_0+A(z_1,z_2,z_3)^T$, with $A=(\phi(i_2),\phi(i_3),\phi(i_5))$. Since $\phi$ is invertible, $A$ is in $\Gl(3,\C)$. Moreover the inner product $\langle,\rangle$ on $\OO$ induces the standard real inner product on $\R^6\cong\C^3$ which is preserved by $\phi$, hence $A$ is in $\U(3)=\SO(6)\cap \Gl(3,\C)$.

Next we claim that two vectors $v, w \in \C^3$ are orthonormal for the standard Hermitian form  if and only if $(i_1, v, w)$ is a Cayley triple.  Indeed, given two octonions $v$ and $w$ are orthonormal for the Hermitian form if and only if $(v, i_1 v, w, i_1 w)$ are orthonormal for the standard real form.  Since they are orthogonal to $i_1$ itself, the result follows. Since $\eG$ acts simply transitively on Cayley triples, ${\eG}_{i_1}$ acts simply transitively on Cayley triples of the form $(i_1, v, w)$. We conclude that ${\eG}_{i_1}$ acts simply transitively on the pairs of orthonormal vectors in $\C^3$. 

We claim that ${\eG}_{i_1} \subseteq \SU(3)$; to do this, it remains to show that every element $A$ as above has determinant one.  Since $A \in \U(3)$, by the spectral theorem there exists an orthonormal eigenbasis for $A$.  By the preceding paragraph, there is an element $\psi \in {\eG}_{i_1}$ such that $B:=\psi^{-1} A \psi$ has $i_2, i_3$ as eigenvectors: let the eigenvalues be $\lambda_1, \lambda_2$ in $\mathbb(O)$.  Now we compute that $B(i_5) = B(i_2 i_3) = B(i_2) B(i_3)=(\lambda_1 i_2)(\lambda_2 i_3) = \overline{\lambda_1 \lambda_2} i_5$.  Thus the determinant of $A$ is one, as required.

Since ${\eG}_{i_1} \subseteq \SU(3)$ and ${\eG}_{i_1}$ acts transitively on the pairs of orthonormal vectors in $\mathbb{C}^3$, it follows that ${\eG}_{i_1} = \SU(3)$, as a transformation in $\SU(3)$ is determined by the image of the first two basis vectors of $\C^3$. As a result the real isotropy representation is just the standard 6-dimensional representation $\mu_3^{\mathbb{R}}$ of $SU(3)$.\\

\paragraph{\textbf{The action of $\Spin(7)$ on $S^7$}:} The spin group $\Spin(n)$ is the connected double covering of $\SO(n)$. In order to understand the action of $\Spin(7)$ and $\Spin(9)$ on spheres, we need to study the real spin representation of these groups.  Explicitly, $\Spin(n)$ can be realised in the following way. Let $\{e_i\}_{i=1}^{n}$ be an orthonormal basis of $(\R^n,\langle\cdot,\cdot\rangle)$ (with the usual scalar product). Recall that the Clifford algebra $Cl_n$ is the algebra over $\R$ generated by the vectors $\{e_i\}_{i=1}^{n}$ satisfying the relations 
$e_i\cdot e_j+e_j\cdot e_i=-2\delta_{ij}.$ Then $\Spin(n)$ is the subgroup of $Cl_n$ given by
$$\Spin(n):=\{x_1\cdot x_2\dots x_{2k}:x_i\in\R^n, \|x_i\|=1, i\in(1,\dots 2k)\}\subset Cl_n$$
Note that the double covering is given by $\lambda: \Spin(n)\rightarrow \SO(n),\, s\mapsto \lambda(s), \lambda(s)x=sxs^{-1}.$ 

The spin representation $\Delta_n:\Spin(n)\rightarrow \Gl(\Sigma_n)$  is obtained by restricting an irreducible representation of $Cl_n$ to $\Spin(n)$. In the case $n=1\mod 4$ (and hence for $n=9$), the Clifford algebra $Cl_{4k+1}$ has a single irreducible representation, which splits under $\Spin(4k+1)$ into two equivalent representations, either one of which defines $\Sigma_n$. 

In the case $n=3 \mod 4$ (and hence for $n=7$),  the representation of the Clifford algebra $Cl_{4k+3}$ is a sum of two inequivalent irreducible representations. The restriction to $\Spin(4k+3)$ is a unique representation, independent of which irreducible representation of $Cl_{4k+3}$ is used.  In particular, in dimension $n=7\mod 8$, $Cl_{8k+7}=\R(2^{3+4k})\oplus\mathbb{R}(2^{3+4k})$ and $\Sigma_n=\R^{2^{3+4k}}$. From this we get immediately the action of $\Spin(7)$ on $\R^8\cong \OO$ and hence on the unit sphere  $S^7\cong \{v\in\OO\cong \R^8:\|v\|=1\}$.

We claim that the action of $\Spin(7)$ on $S^7$ is transitive with isotropy group $\eG$.  Here, $\eG$ is realised inside of $\Spin(7)$ as follows: $\eG < \SO(7)$ as above and since $\eG$ is simply connected, this inclusion lifts to $\eG < \Spin(7)$.  

To prove the claim, we first show that $\Spin(7)$ acts transitively on $S^7$. To see this we realise the spin representation explicitly in the following known way. There is an action of $Cl(\R^7)\cong Cl(\textrm{Im}(\mathbb{O}))$ on $\OO$ given by the left multiplication $\textrm{Im}(\mathbb{O})\times \OO \to \OO$ of imaginary octonions. The fact that this defines an action of $Cl(\R^7)$ follows from explicitly checking the Clifford relations---a consequence of the multiplication table of $\OO$.  As this is a nontrivial module, it must be one of the irreducible modules of dimension eight of $Cl(\R^7)$, therefore its restriction to $\Spin(7)$ is the spin representation.  The Lie algebra $\mathfrak{so}(7) \subseteq Cl(\R^7)$ is the span of commutators of the generators under Clifford multiplication. Let us compute $\mathfrak{so}(7) \cdot 1$, where $1 \in \OO$ is considered a unit vector, hence in $S^7$.  For $v, w \in \R^7$, so $v \wedge w \in \mathfrak{so}(7)$, we have $(v \wedge w) \cdot 1 = vw - wv$, the commutator in the octonions.  We see from the multiplication table that commutators of imaginary octonions span the imaginary octonions again.   Thus, $\mathfrak{so}(7) \cdot 1 = \R^7$, the imaginary octonions.  Therefore, the orbit of $1$ under $\Spin(7)$ has dimension $7$. Since $\Spin(7)$ is compact, the orbit is a closed submanifold of $S^7$, hence all of $S^7$. This implies transitivity.

We now compute the isotropy. Note that $\eG < \Spin(7)$ acts on the spin representation $\R^8$.  Since the only nontrivial irreducible representation of $\eG$ of dimension $\leq 8$ has dimension $7$, and all finite-dimensional representations of $\eG$ are completely reducible,  it follows that there is a trivial $\eG$-subrepresentation of $\R^8$. That means there is a fixed vector, hence a unit fixed vector $v \in S^7$. Let $H$ be the isotropy group at $v$ then we have $\eG < H$ by definition of $v$. Since $\Spin(7)$ acts transitively on $S^7$,  we have $\dim H = 21-7=14 = \dim \eG$.  Thus, $\eG < H$ is an open subgroup, and as $G_2$ is compact, it is the connected component of the identity.  Finally, since $\Spin(7)/H \cong S^7$ is simply connected, by the long exact sequence on homotopy groups, $|\pi_0(H)|=1$, so $H=\eG$ as desired.

Since the isotropy representation is a nontrivial seven-dimensional representation of $G_2$, it has to be isomorphic to $\zeta$. \\

\paragraph{\textbf{The action of $\Spin(9)$ on $S^{15}$.}} In dimension $n=1\mod 8$, we have $Cl_{8k+1}=\C(2^{4k})$, which after restricting to $\Spin(8k+1)$ yields $\Sigma_n=\R^{2^{4k}}$. Therefore in dimension 9, $\Sigma_9=\R^{16}$ and $\Spin(9)$ acts via the spin representation on the unit sphere $S^{15}\subset\R^{16}\cong \mathbb{O}^2$. Moreover, we identify: \begin{equation*}
\mathbb{R}^9\cong \mathbb{R}\oplus\mathbb{O}=\Biggl\{i\begin{pmatrix} r & L_z\\L_{\bar{z}} & -r\end{pmatrix} : r\in\mathbb{R}, z\in\mathbb{O}\Biggr\} \subseteq \textrm{End}_{\C}(\C \otimes \OO^2),
\end{equation*}where $L_z$ denotes the map associated to the left multiplication by $z$ in $\OO$. There is, similarly to $\Spin(7)$, an action of $Cl(\R\oplus\OO)$ on $\C \otimes \OO^2$ given by left multiplication.
Note that in this case, $\mathbb{C} \otimes \mathbb{O}^2$ is the unique irreducible representation of the real Clifford algebra $Cl(\mathbb{R}\oplus\mathbb{O})$ and the action furnishes an isomorphism $Cl(\mathbb{R}\oplus\mathbb{O}) \cong \textrm{End}_{\mathbb{C}}(\C \otimes \OO^2)$.  Restricting this to the even Clifford algebra, it splits into two isomorphic real representations, $\OO^2$ and $i \OO^2$. Hence, these both identify with $\Sigma_9$.

We can compute directly, as in the case of $\Spin(7)$, that  $\Spin(9)$ acts transitively on $S^{15}$.  Take the vector $(1,0) \in \OO^2$. Since the orbit $\Spin(9) \cdot (1,0) \subseteq S^{15}$ is closed, to show it is all of the sphere it suffices to prove that the tangent space $\mathfrak{so}(9) \cdot (1,0)$ has dimension at least $15$.  Now, $\mathfrak{so}(9)$ is the span in $Cl(\R \oplus \OO)$ of commutators of generators.  We compute that
\begin{equation*}
\Biggl[ i\begin{pmatrix} r & L_u\\L_{\bar{u}} & -r\end{pmatrix}, 
i\begin{pmatrix} r' & L_v\\L_{\bar{v}} & -r'\end{pmatrix} \Biggr]
= 
-\begin{pmatrix} L_u L_{\bar v} - L_v L_{\bar u} & 2(rL_v - r' L_u) \\
2(r' L_{\bar u} - rL_{\bar v}) & L_{\bar u} L_v - L_{\bar v} L_{u}
\end{pmatrix}.
\end{equation*}
Hence the tangent space is the span of the vectors
\begin{equation}\label{spin9-tgt}
((r,u) \wedge (r',v)) \cdot (1,0) = -(u\bar v - v \bar u, 2(r'\bar u - r \bar v)), r,r' \in \R, u,v \in \OO.
\end{equation}
Taking $r=1$ and $r'=0, u=0$, we obtain $(0,-2\overline{v})$: this spans $(0,\OO)$.
Taking $u=1_{\OO}, r=r'=0$, we obtain $(\bar v - v, 0)$: this spans $(\text{Im}(\OO),0)$.  So the overall span is at least $\text{Im}(\OO) \oplus \OO$, which is the tangent space to $S^{15}$ and we are done.

Note that we can explicitly compute the isotropy Lie algebra, call it $\mathfrak{h}$, of $(1,0)$ in $\mathfrak{so}(9)$.  We claim that this is:
\begin{equation}\label{spin9-iso}
\text{span}\bigl(T_{u,v} := u \wedge v + \frac{1}{2}([u,v] \wedge 1_{\OO}) \mid u,v \in \text{Im}(\OO) \bigr) \subseteq \wedge^2 \OO  = \mathfrak{so}(8) \subseteq \wedge^2(\R \oplus \OO)=\mathfrak{so}(9).
\end{equation}
Indeed, plugging $T_{u,v}$ into \eqref{spin9-tgt} with $r=r'=0$ we get zero, since $(uv-vu) -\frac{1}{2}([u,v]+[u,v]) = 0$.  So $T_{u,v} \in \mathfrak{h}$.  On the other hand, these elements span a $21=(36-15)$-dimensional subspace, which is the correct dimension. More explicitly, a complementary subspace to \eqref{spin9-iso} is $\R \wedge \OO \oplus \mathrm{Re}(\OO) \wedge \mathrm{Im}(\OO)$, which we already showed acts injectively on $(1,0)$, with image the tangent space to $S^{15}$.

The isotropy representation is the tangent space $\text{Im}(\OO) \oplus \OO$. The action of the isotropy Lie algebra decomposes into two subrepresentations, $(\text{Im}(\OO),0) \oplus (0,\OO)$.  The action on the first factor induces a homomorphism $\varphi: \mathfrak{h} \to \mathfrak{gl}(7)$. Explicitly, this action is given on $\text{Im}(\OO)$ as follows, for $u,v,w \in \text{Im}(\OO)$:
\begin{equation}\label{action-so7}
T_{u,v} \cdot w = u \cdot (v \cdot w) - v \cdot (u \cdot w) - ([u,v]) \cdot w,
\end{equation}
with multiplication taken in the octonions.  This can be interpreted as ``the failure of left multiplication of $\mathrm{Im}(\OO)$ on itself to define a Lie algebra representation'' (although note that $\mathrm{Im}(\OO)$ is not a Lie algebra under commutators).   

We claim that the image of this homomorphism lies in $\mathfrak{so}(7)$.  Suppose $u=i_1, v=i_2, w = i_j$. Then $u(vw) = \pm (uv)w=\mp (vu)w = -v(uw)$, with the sign $\pm$ positive if $j \in \{1,2,4\}$, and  negative otherwise.  It follows that the action of $T_{u,v}$ is given by $i_3 \mapsto -i_6 \mapsto i_3$ and $i_5 \mapsto i_7 \mapsto -i_5$.  In other words, $\varphi(T_{u,v})$ is zero on the three-dimensional space spanned by $i_1, i_2, i_4$ (a line in the octonion Fano plane), and acts by $90^\circ$ rotation on the two two-dimensional spaces spanned by $i_3,i_6$ and by $i_5,i_7$ (which with $i_4$ form the other three lines through $i_4$).  In particular, $\varphi(T_{u,v}) \in \mathfrak{so}(7)$ (it is skew-symmetric). Through automorphisms of the Fano plane, this describes $\varphi(T_{i_k, i_\ell})$ for all indices $k,\ell$. Indeed, each transformation is zero on the span of $i_k, i_\ell$, and $i_k i_\ell$ (which represent a line in the octonionic Fano plane), and is a $90^\circ$ rotation on each of the complementary two-dimensional spaces of the form $\text{Span}(i_m, (i_k i_\ell) i_m)$ (representing the two other lines in the octonionic Fano plane containing $\pm i_k i_\ell$). 

We claim furthermore that this gives an isomorphism $\mathfrak{h} \to \mathfrak{so}(7)$. Since the source and target have the same dimension, it suffices to show that the $\varphi(T_{u,v})$ are linearly independent: since $T_{i_k, i_\ell}$ acts on some of the $i_j$ as multiplication by
$i_k i_\ell$ and by zero on others, it suffices to check that the three choices of indices $(k,\ell)$ with the same product $i_k i_\ell$ (i.e., the representing the three lines in the octonionic Fano plane through $i_k i_\ell$) give linearly independent elements $\varphi(T_{i_k, i_\ell})$. That is, $\varphi(T_{i_1, i_2}), \varphi(T_{i_6, i_3})$, and $\varphi(T_{i_5, i_7})$ should be linearly independent (and cyclic permutations of $(1,2,3,4,5,6,7)$).  This is true: if $\lambda_1 T_{i_1, i_2} + \lambda_2 T_{i_6, i_3} + \lambda_3 T_{i_5, i_7}$ acts by zero on $i_1, i_2$, then $\lambda_2=-\lambda_3$; similarly to act by zero on $i_6, i_3$ we need $\lambda_3=-\lambda_1$, and finally to act by zero on $i_5,i_7$ we need $\lambda_3=-\lambda_1$. So $\lambda_1=\lambda_2=\lambda_3=0$. Similarly, applying cyclic permutations of the indices $(1,2,3,4,5,6,7)$ we get that a linear combination of the $T_{i_k,i_\ell}$ which acts by zero is the zero combination.

As a result, the isotropy representation $\text{Im}(\OO) \oplus \OO$ of $\mathfrak{h}$ realises $\mathfrak{h} \cong \mathfrak{so}(7)$ with $\mathrm{Im}(\OO)$ the standard $7$-dimensional representation.  Note that the isotropy Lie subgroup $H$ must be simply-connected by the long exact sequence on homotopy groups for $H \to \Spin(9) \to S^{15}$, hence $H \cong \Spin(7)$.  Therefore the other summand $(0,\OO)$ of the isotropy representation must be isomorphic to the spin representation, otherwise the isotropy representation would not be faithful. 

This is a contradiction, because $H$ acts effectively on $S^{15}$.  More generally, if $G$ acts effectively on a compact connected manifold $M$ and $H$ is the stabilizer of $p\in M$ is compact, we claim that the isotropy representation is faithful.  Let $K<H$ be the kernel of the isotropy representation. 
Then $(T_p M)^K = T_p (M^K)$, since there is an $K$-equivariant isomorphism of a neighbourhood of $p$ in $M$ and of $p$ in $T_p M$  (or because we can take an $K$-invariant metric on $M$ and then $K$ preserves the exponential flow beginning at $p$, hence all geodesics on $X$ through $p$).  As a result, since $K$ acts trivially on $T_p M$ and $M$ is compact and connected, $M^K$ is an open and closed subset of $M$, hence all of $M$.  So $K$ acts trivially on $M$. As the action of $G$ is effective, $K$ is trivial, as desired.

Thus, this explicitly realizes $H \cong \Spin(7)$ and the isotropy representation as a direct sum of the standard and spin representations.


\end{document}